\documentclass[a4paper,12pt]{article}
\usepackage{amssymb,fullpage,amsmath,amsthm}

\newtheorem{thm}{Theorem}[section]

\newtheorem{cor}[thm]{Corollary}

\newtheorem{lem}[thm]{Lemma}
\newtheorem{propn}[thm]{Proposition}
\newtheorem{rem}{Remark}

\begin{document}

\title{Random walks on Galton-Watson trees with infinite variance offspring distribution conditioned to survive}
\author{David Croydon\footnote{Dept of Statistics,
University of Warwick, Coventry CV4 7AL, UK;
{d.a.croydon@warwick.ac.uk}.} and Takashi Kumagai\footnote{Dept of Mathematics, Kyoto University, Kyoto 606-8502, Japan; {kumagai@math.kyoto-u.ac.jp}.}}
\date{4 August 2008}
\maketitle

\begin{abstract}
We establish a variety of properties of the discrete time simple random walk on a Galton-Watson tree conditioned to survive when the offspring distribution, $Z$ say, is in the domain of attraction of a stable law with index $\alpha\in(1,2]$. In particular, we are able to prove a quenched version of the result that the spectral dimension of the random walk is $2\alpha/(2\alpha-1)$. Furthermore, we demonstrate that when $\alpha\in(1,2)$ there are logarithmic fluctuations in the quenched transition density of the simple random walk, which contrasts with the log-logarithmic fluctuations seen when $\alpha=2$. In the course of our arguments, we obtain tail bounds for the distribution of the $n$th generation size of a Galton-Watson branching process with offspring distribution $Z$ conditioned to survive, as well as tail bounds for the distribution of the total number of individuals born up to the $n$th generation,  that are uniform in $n$.
\end{abstract}

\section{Introduction}

This article contains an investigation of simple random walks on Galton-Watson trees conditioned to survive, and we will start by introducing some relevant branching process and random walk notation. Let $Z$ be a critical ($\mathbf{E}Z=1$) offspring distribution in the domain of attraction of a stable law with index $\alpha\in (1,2]$, by which we mean that there exists a sequence $a_n\uparrow \infty$ such that
\begin{equation}\label{attraction}
\frac{Z[n]-n}{a_n}\buildrel{d}\over{\rightarrow} X,
\end{equation}
where $Z[n]$ is the sum of $n$ i.i.d. copies of $Z$ and $\mathbf{E}(e^{-\lambda X})=e^{-\lambda^\alpha}$. Note that, by results of \cite{Feller2}, Chapters XIII and XVII, this is equivalent to $Z$ satisfying
\begin{equation}\label{genfun}
\mathbf{E}(s^Z)=s+(1-s)^{\alpha}L(1-s),\hspace{20pt}\forall s\in (0,1),
\end{equation}
where $L(x)$ is slowly varying as $x\rightarrow 0^+$, and the non-triviality condition $\mathbf{P}(Z=1)\neq 1$ holding. Denote by $(Z_n)_{n\geq 0}$ the corresponding Galton-Watson process, started from $Z_0=1$. It has been established (\cite{Slack}, Lemma 2) that if $p_n:=\mathbf{P}(Z_n>0)$, then
\begin{equation}\label{probdecay}
p_n^{\alpha-1}L(p_n)\sim \frac{1}{(\alpha-1)n},
\end{equation}
as $n\rightarrow\infty$, where $L$ is the function appearing in (\ref{genfun}). It is also well known that the branching process $(Z_n)_{n\geq 0}$ can be obtained as the generation size process of a Galton-Watson tree, $T$ say, with offspring distribution $Z$. In particular, to construct the random rooted graph tree $T$, start with an single ancestor (or root), and then suppose that individuals in a given generation have offspring independently of the past and each other according to the distribution of $Z$, see \cite{rrt}, Section 3, for details. The vertex set of $T$ is the entire collection of individuals, edges are the parent-offspring bonds, and $Z_n$ is the number of individuals in the $n$th generation of $T$. By (\ref{probdecay}), it is clear that $T$ will be a finite graph $\mathbf{P}$-a.s. However, in \cite{Kesten}, Kesten showed that it is possible to make sense of conditioning $T$ to survive or ``grow to infinity''. More specifically, there exists a unique (in law) random infinite rooted locally-finite graph tree $T^*$ that satisfies, for any $n\in\mathbb{Z}_+$,
\begin{equation}\label{tstar}
\mathbf{E}\left(\phi(T^*|_n)\right)=\lim_{m\rightarrow\infty}\mathbf{E}\left(\phi(T|_n)|Z_{m+n}>0\right),
\end{equation}
where $\phi$ is a bounded function on finite rooted graph trees of $n$ generations, and $T|_n$, $T^*|_n$ are the first $n$ generations of $T$, $T^*$ respectively. It is known that, for any $n\in\mathbb{Z}_+$, we can also write
\begin{equation}\label{sizebias}
\mathbf{E}\left(\phi(T^*|_n)\right)=\mathbf{E}\left(\phi(T|_n)Z_n\right),
\end{equation}
for any bounded function $\phi$ on finite rooted graph trees of $n$ generations (see \cite{Kesten}, Lemma 1.14, for example), which demonstrates that the the law of the first $n$ generations of $T^*$ is simply a size-biased version of the law of the first $n$ generations of the unconditioned tree $T$. Finally, from the characterisation of $T^*$ at (\ref{tstar}), it is clear that the generation size process $(Z^*_n)_{n\geq 0}$ of $T^*$ is precisely the $Q$-process associated with $(Z_n)_{n\geq 0}$ (see \cite{A-N}, Section I.14), which is commonly referred to as the Galton-Watson process conditioned to survive.

Given a particular realisation of $T^*$, let $X=((X_m)_{m\geq 0}, P_x^{T^*},x\in T^*)$ be the discrete time simple random walk on $T^*$. Define a measure $\mu^{T^*}$ on $T^*$ by setting $\mu^{T^*}(A)=\sum_{x\in A}\mathrm{deg}_{T^*}(x)$, where $\mathrm{deg}_{T^*}(x)$ is the graph degree of the vertex $x$ in $T^*$. The measure $\mu^{T^*}$ is invariant for $X$, and the transition density of $X$ with respect to $\mu^{T^*}$ is given by
\[p_m^{T^*}(x,y):=\frac{{P}^{T^*}_x(X_m=y)}{\mu^{T^*}(\{y\})},\hspace{20pt}\forall x,y\in T^*,m\in \mathbb{Z}_+.\]
Throughout this article, we use the notation
\[\tau_R:=\min\{m:d_{T^*}(X_0,X_m)=R\},\]
where $d_{T^*}$ is the usual shortest-path graph distance on $T^*$, to represent the first time that $X$ has traveled a distance $R\in \mathbb{N}$ from its starting point.

The behaviour of $X$, started from the root $\rho$ of $T^*$, and $(\tau_R)_{R\geq 1}$ was first considered in \cite{Kesten}, where Kesten showed that, under the annealed law
\[\mathbb{P}:=\int P^{T^*}_\rho(\cdot) d\mathbf{P},\]
if the offspring distribution has finite variance, then the rescaled height-process defined by $(n^{-1/3}d_{T^*}(\rho,X_{\lfloor nt\rfloor}))_{t\geq 0}$ converges weakly as $n\rightarrow\infty$ to a non-trivial continuous process (the full proof of this result appeared in \cite{Kestenunpub}). In \cite{Kesten}, it was also noted that
\begin{equation}\label{kest}
\lim_{\lambda\rightarrow\infty}\inf_{R\geq1}\mathbb{P}\left(\lambda^{-1}R^{\frac{2\alpha-1}{\alpha-1}}\leq \tau_R\leq \lambda R^{\frac{2\alpha-1}{\alpha-1}}\right)=1,
\end{equation}
whenever the offspring distribution is in the domain of {\it normal} attraction of a stable law with index $\alpha\in(1,2]$, by which we mean that there exists a constant $c\in(0,\infty)$ such that (\ref{attraction}) occurs with $a_n=cn^{\frac{1}{\alpha}}$ (the full proof was given in the case $\alpha=2$ only). More recently, in the special case when $Z$ is a binomial random variable, a detailed investigation of $X$ was undertaken in \cite{BarKum}, where a variety of bounds describing the quenched (almost-sure) and expected behaviour of the transition density and displacement, $d_{T^*}(\rho,X_m)$, were established. Many of these results have since been extended to the general finite variance offspring distribution case, see \cite{KF}.

In this article, we continue the above work by proving distributional, annealed and quenched bounds for the exit times, transition density and displacement of the random walk on $T^*$ for general offspring distributions satisfying (\ref{attraction}). Similarly to the arguments of \cite{BarKum} and \cite{KF}, to deduce properties of the random walk, it will be important to estimate geometric properties of the graph $T^*$ such as the volume growth and resistance across annuli (when $T^*$ is considered as an electrical network with a unit resistor placed on each edge). In particular, once we have established suitable volume growth and resistance bounds, we can apply the results proved for general random graphs in \cite{KM} to obtain many of our random walk conclusions. It should be noted that the techniques applied in \cite{KM} are simple extensions of those developed in \cite{BJKS}, which apply in our case when $\alpha=2$.

In terms of the branching process, we are required to derive suitably sharp estimates on the tails of the distributions of $Z_n^*$ and $\sum_{m=0}^nZ_m^*$, which we do in Section \ref{bpp}. In \cite{KF}, bounds for these quantities were obtained in the finite variance offspring distribution case using moment estimates
which fail in the more general case that we consider here. We are able to overcome this problem using more technical arguments, which involve analysis of the generating functions of the relevant random variables. The statement of our results depends on the non-extinction probabilities of the branching process $(p_n)_{n\geq 0}$ via a ``volume growth function'' (see Section 3 for a justification of this title), $v:\mathbb{R}_+\rightarrow\mathbb{R}_+$, which is defined to satisfy $v(0)=0$,
\begin{equation}\label{vgf}
v(R):=Rp_R^{-1},\hspace{20pt}\forall R\in\mathbb{N},
\end{equation}
and is linear in-between integers. Our first result yields the tightness of the suitably rescaled distributions of $(\tau_R)_{R\geq 1}$, $(E^{T^*}_\rho\tau_R)_{R\geq 1}$, $(p^{T^*}_{2m}(\rho,\rho))_{m\geq 1}$ and $(d_{T^*}(\rho,X_m))_{m\geq 1}$ with respect to the appropriate measures; along with all the subsequent theorems of this section, it is proved in Section \ref{proofs}.

{\thm \label{probthm} The random walk on $T^*$ satisfies
\begin{equation}\label{kestens}
\lim_{\lambda\rightarrow\infty}\inf_{R\in\mathbb{N}}\mathbb{P}\left(\lambda^{-1}h(R)\leq \tau_R\leq \lambda h(R)\right)=1,
\end{equation}
\[\lim_{\lambda\rightarrow\infty}\inf_{R\in\mathbb{N}}\mathbf{P}\left(\lambda^{-1}h(R)\leq E^{T^*}_\rho\tau_R\leq \lambda h(R)\right)=1,\]
\[\lim_{\lambda\rightarrow\infty}\inf_{m\in\mathbb{N}}\mathbf{P}\left(\lambda^{-1}\leq v(\mathcal{I}((m))p^{T^*}_{2m}(\rho,\rho)\leq \lambda \right)=1,\]
\begin{equation}\label{dlim}
\lim_{\lambda\rightarrow\infty}\inf_{m\in\mathbb{N}}\mathbb{P}\left(\lambda^{-1} \mathcal{I}(m)\leq1+d_{T^*}(\rho,X_m)\leq \lambda \mathcal{I}(m)\right)=1,
\end{equation}
where $h(R):=Rv(R)$ and $\mathcal{I}(m):=h^{-1}(m)$.}
\bigskip

We remark that, from (\ref{probdecay}), we have that $v(R)=R^{\frac{\alpha}{\alpha-1}}\ell(R)$ for some function $\ell$ which is slowly varying as $R\rightarrow\infty$ (see Lemma \ref{pnprop} below). Thus the functions bounding the exit time, transition density and displacement in the above result are of the form:
\[h(R)=R^{\frac{2\alpha-1}{\alpha-1}}\ell_1(R),\hspace{10pt}v(\mathcal{I}(m))=m^{\frac{\alpha}{2\alpha-1}}\ell_2(m),\hspace{10pt}\mathcal{I}(m)=m^\frac{\alpha-1}{2\alpha-1}\ell_3(m),\]
where $\ell_1$, $\ell_2$ and $\ell_3$ are slowly varying at infinity. In particular, when $Z$ is in the domain of normal attraction of a stable law with index $\alpha$, then we have that $p_n\sim cn^{-\frac{1}{\alpha-1}}$ for some constant $c$, and hence (\ref{kestens}) provides an alternative proof of the result of Kesten's stated at (\ref{kest}). We highlight the fact that the $\alpha$ of Kesten's article corresponds to our $\alpha-1$.

The annealed bounds that we are able to obtain include the following. Further off-diagonal estimates for the transition density, which extend the estimate at (\ref{anntd}) are presented in Section \ref{offd}.

{\thm \label{annthm} For every $\beta\in(0,\alpha-1)$, $\gamma\in(0,1-\alpha^{-1})$ and $\delta\in(0,\alpha^{-1})$, there exist constants $c_1,\dots,c_6\in(0,\infty)$ such that
\begin{equation}\label{anntau}
c_1h(R)^\beta\leq\mathbf{E}\left(\left(E^{T^*}_\rho\tau_R\right)^\beta\right)\leq c_2h(R)^\beta,\hspace{20pt}\forall R\in\mathbb{N},
\end{equation}
\begin{equation}\label{anntd}
c_3v(\mathcal{I}(m))^{-\gamma}\leq\mathbf{E}\left(p_{2m}^{T^*}(\rho,\rho)^\gamma\right)\leq c_4v(\mathcal{I}(m))^{-\gamma},\hspace{20pt}\forall m\in\mathbb{N},
\end{equation}
\begin{equation}\label{annd}
c_5\mathcal{I}(m)^\delta\leq \mathbb{E}\left(d_{T^*}(\rho,X_m)^\delta\right)\leq \mathbb{E}\left(\max_{0\leq k\leq m} d_{T^*}(\rho,X_k)^\delta\right)\leq c_6\mathcal{I}(m)^\delta,\hspace{20pt}\forall m\in\mathbb{N}.
\end{equation}}
\bigskip

In the finite variance case, it is known that (\ref{anntau}) and (\ref{anntd}) hold with $\beta,\gamma=1$ (see \cite{KF}, Theorem 1.1). Furthermore, in \cite{BarKum}, it was established that when the offspring distribution is binomial, then (\ref{annd}) holds with $\delta=1$. The proofs of (\ref{anntau}) and the corresponding results in \cite{BarKum} and \cite{KF} all rely on the bound $E^{T^*}_\rho\tau_R\leq 2(R+1)\sum_{m=0}^{R+1}Z_m^*$. However, the right-hand side of this expression has infinite expectation under $\mathbf{P}$ when $\alpha\in(1,2)$, and so we can not use the same technique to deduce the result for $\beta=1$ in general. A similar problem occurs in the proof of (\ref{anntd}), where, to establish the result for $\gamma=1$, we need an estimate on the negative moments of $\sum_{m=0}^{R}Z_m^*$ of orders larger than we can prove. We cannot prove if (\ref{anntau}) and (\ref{anntd}) actually fail to hold or not in general when $\beta=\gamma=1$. We also do not know whether, when $\delta=1$, the expectations in (\ref{annd}) can be bounded above by a multiple of $\mathcal{I}(m)$ uniformly in $m$ in general.

In addition to the above annealed bounds, we will also establish quenched bounds for the random walk on $T^*$ as follows. Note that part (b) implies that for $\mathbf{P}$-a.e. realisation of $T^*$, the random walk on $T^*$ is recurrent.

{\thm\label{qthm} There exist constants $a_1,\dots,a_4\in(0,\infty)$ such that for $\mathbf{P}$-a.e. realisation of $T^*$ the following properties are satisfied.\\
(a) If $x\in T^*$, then $P^{T^*}_x$-a.s.,
\[h(R)(\log R)^{-a_1}\leq \tau_R\leq h(R)(\log R)^{a_1},\hspace{20pt}\mbox{for large $R$.}\]
Moreover,
\[h(R)(\log R)^{-a_2}\leq E_x^{T^*}\tau_R\leq h(R)(\log R)^{a_2},\hspace{20pt}\mbox{for large $R$.}\]
(b) If $x\in T^*$, then
\[v(\mathcal{I}(m))^{-1}(\log m)^{-a_3}\leq p_{2m}^{T^*}(x,x)\leq v(\mathcal{I}(m))^{-1}(\log m)^{a_3},\hspace{20pt}\mbox{for large $m$.}\]
(c) If $x\in T^*$, then $P^{T^*}_x$-a.s.,
\[\mathcal{I}(m)(\log m)^{-a_4}\leq \max_{0\leq k\leq m} d_{T^*}(x,X_k)\leq \mathcal{I}(m)(\log m)^{a_4},\hspace{20pt}\mbox{for large $m$.}\]}

From the preceding theorem we are easily able to calculate the exponents of the leading order polynomial terms governing the exit time, transition density decay and maximum displacement. We are also able to determine the exponent according to which the size of the range of the simple random walk grows.

{\thm \label{exponents} For $\mathbf{P}$-a.e. realisation of $T^*$ , we have that
\begin{equation}\label{exptau}
\lim_{R\rightarrow\infty}\frac{\log \tau_R}{\log R}=\lim_{R\rightarrow\infty}\frac{\log E^{T^*}_x(\tau_R)}{\log R}=\frac{2\alpha-1}{\alpha-1},
\hspace{20pt}P^{T^*}_x\mbox{-a.s. for every $x\in T^*$},
\end{equation}
\begin{equation}\label{exptd}
d_s(T^*):=\lim_{m\rightarrow\infty}\frac{-2\log p_{2m}^{T^*}(\rho,\rho)}{\log m}=\frac{2\alpha}{2\alpha-1},
\end{equation}
\begin{equation}\label{expd}
\lim_{m\rightarrow\infty}\frac{\log \max_{0\leq k\leq m} d_{T^*}(x,X_k)}{\log m}=\frac{\alpha-1}{2\alpha-1},
\hspace{20pt}P^{T^*}_x\mbox{-a.s. for every $x\in T^*$},
\end{equation}
and if the range $W=(W_m)_{m\geq 0}$ of the simple random walk is defined by setting $W_m:=\{X_0,\dots,X_m\}$, then
\[\lim_{m\rightarrow\infty}\frac{\log \mu^{T^*}(W_m)}{\log m}=\lim_{m\rightarrow\infty}\frac{\log \# W_m}{\log m}=\frac{\alpha}{2\alpha-1}, \hspace{20pt}P^{T^*}_x\mbox{-a.s. for every $x\in T^*$}.\]}

We remark that the quantity $d_s(T^*)$ introduced in the above result is often taken as a definition of the spectral dimension of (the random walk on) $T^*$. Famously, in \cite{AO}, Alexander and Orbach conjectured that the random walk on the ``infinite cluster [of the bond percolation model on $\mathbb{Z}^d$] at the critical percolation concentration'' has spectral dimension $4/3$, independent of the spatial dimension $d$. Although a precise definition of the percolative graph considered did not appear in \cite{AO}, it is now commonly interpreted as the incipient infinite cluster of critical percolation, as constructed in \cite{Kesten} for $d=2$ and in \cite{HJ} for large $d$. The validity of the Alexander-Orbach conjecture has since been challenged in small dimensions \cite{Hughes}, but it is still widely believed to hold above a certain critical dimension (it has in fact been established in the case of spread-out oriented percolation in high dimensions \cite{BJKS}). Justification for such a conviction is provided by the known results about the geometry of the incipient infinite percolation cluster on $\mathbb{Z}^d$ (see \cite{HHS} and \cite{HJ}, for example), which suggest that it is closely related to the incipient infinite percolation cluster on an $N$-ary tree or, equivalently, the Galton-Watson tree with binomial offspring distribution, parameters $N$ and $p=N^{-1}$, conditioned to survive, where versions of the Alexander-Orbach conjecture are known to hold. For example, for the more general offspring distributions considered in \cite{Kesten}, Kesten explains how the result at (\ref{kest}) implies that the Alexander-Orbach conjecture holds for $T^*$ if and only if $\alpha=2$, presenting the discussion in terms of distributional scaling exponents. Contributing to these developments, it is worthwhile to observe that the limit result at (\ref{exptd}) yields a quenched version of this dichotomy between the cases $\alpha=2$, where the Alexander-Orbach conjecture holds, and $\alpha\in(1,2)$, where it does not.

Finally, in Section \ref{fluct}, we investigate the fluctuations in the volume growth and the quenched transition density of the simple random walk on $T^*$. In particular, when $\alpha\in(1,2)$ we show that the volume of a ball of radius $R$, centered at $\rho$, has logarithmic fluctuations about the function $v(R)$, $\mathbf{P}$-a.s., and when $\alpha=2$ there are log-logarithmic fluctuations, $\mathbf{P}$-a.s. It follows from estimates in Section \ref{bpp} and \cite{BarKum} that these results are sharp up to exponents. We also note that these asymptotic results are analogous to the results proved for the related stable trees in \cite{treemeas}, where it is shown that the Hausdorff measure of a stable tree with index $\alpha\in(1,2)$ has logarithmic corrections, in contrast to the log-logarithmic corrections seen when $\alpha=2$. Furthermore, by standard arguments, it follows that, with positive probability, the quenched transition density of the simple random walk on $T^*$ has logarithmic fluctuations when $\alpha\in(1,2)$, and log-logarithmic fluctuations when $\alpha=2$. In general,  these results are also sharp up to exponents in the fluctuation terms.

\section{Branching process properties}\label{bpp}

To establish the properties of the simple random walk on $T^*$ that are stated in the introduction we start by studying the associated generation size process $(Z^*_n)_{n\geq 0}$. That the rescaled sequence $(p_nZ_n^*)_{n\geq 0}$ converges in distribution to a non-zero random variable was proved as \cite{Pakes}, Theorem 4. Furthermore, if we define $(Y^*_n)_{n\geq 0}$ by setting
\[Y^*_n=\sum_{m=0}^nZ^*_m,\]
then it is possible to deduce that $(n^{-1}p_nY^*_n)_{n\geq 0}$ converges in distribution to a non-zero random variable by applying Theorem 1.5 of \cite{Duqinf}, (in fact, \cite{Duqinf}, Theorem 1.5 also provides an alternative description of the limit random variable of $(p_nZ_n^*)_{n\geq 0}$). However, although these results are enough to demonstrate that Theorem \ref{probthm} is true, to deduce the remaining results of the introduction, we need to ascertain tail estimates for $Z^*_n$ and $Y^*_n$ that are uniform in $n$, and that is the aim of this section. We start by stating a moment estimate for the unconditioned Galton-Watson process $(Z_n)_{n\geq 0}$.

{\lem[\cite{FVW}, Lemma 11]\label{moments} For $\beta\in (0,\alpha-1)$, there exists a finite constant $c$ such that
\[\mathbf{E}\left(Z_n^{1+\beta}\right)\leq c p_n^{-\beta},\hspace{20pt}\forall n\in \mathbb{N}.\]}

A polynomial upper bound for the tail of the distribution of $Z_n^*$ near infinity is an easy consequence of this result. When $\alpha\in(1,2)$, we are also able to deduce a polynomial lower bound by first proving a related bound for the generating function of $Z^*_n$.

{\propn\label{zupper} For $\beta_1\in (0,\alpha-1)$, there exists a finite constant $c_1$ such that
\begin{equation}\label{zupperb}
\mathbf{P}\left(Z_n^*\geq \lambda p_n^{-1}\right)\leq c_1 \lambda^{-\beta_1},\hspace{20pt}\forall n\in \mathbb{N},\lambda>0.
\end{equation}
Moreover, for $\alpha\in(1,2)$, $\beta_2>(\alpha-1)/(2-\alpha)$, there exists a strictly positive constant $c_2$ and integer $n_0$ such that
\begin{equation}\label{zupperc}
\mathbf{P}\left(Z_n^*\geq \lambda p_n^{-1}\right)\geq c_2 \lambda^{-\beta_2},\hspace{20pt}\forall n\geq n_0,\lambda\geq 1.
\end{equation}}
\begin{proof} Fix $\beta_1\in (0,\alpha-1)$. Applying the size-biasing result that appears at (\ref{sizebias}), it is possible to deduce that
\[\mathbf{P}\left(Z_n^*\geq \lambda p_n^{-1}\right)=\mathbf{E}\left(Z_n\mathbf{1}_{\{Z_n\geq  \lambda p_n^{-1}\}}\right)\leq \lambda^{-\beta_1}p_n^{\beta_1}\mathbf{E}\left(Z_n^{1+{\beta_1}}\right).\]
Combining this bound and Lemma \ref{moments} yields the upper bound at (\ref{zupperb}).

To prove (\ref{zupperc}), we start by demonstrating that for each $\varepsilon>0$ there exists a constant $c_1$ and integer $n_0$ such that
\begin{equation}\label{gen}
\mathbf{E}\left(e^{-\lambda p_nZ^*_n}\right)\leq 1-c_1\lambda^{\alpha-1+\varepsilon},\hspace{20pt}\forall n\geq n_0,\lambda\in[0,1].
\end{equation}
Let $f(s)=\mathbf{E}(s^Z)$, denote by $f_n$ the $n$-fold composition of $f$, and set
\begin{equation}\label{udef}
U(s)=\lim_{n\rightarrow\infty}\frac{f_n(s)-f_n(0)}{f_n(0)-f_{n-1}(0)},\hspace{20pt}\forall s\in[0,1).
\end{equation}
That such a limit exists for each $s\in[0,1)$ is proved in \cite{Slack}, where it is also established that the resulting function $U$ (which is actually the generating function of the stationary measure of the branching process $Z$) satisfies
\begin{equation}\label{usum}
U(f(s))=U(s)+1,\hspace{20pt}\forall s\in[0,1),
\end{equation}
and we can write
\begin{equation}\label{utail}
U(s)^{-1}=(\alpha-1)(1-s)^{\alpha-1}M(1-s),
\end{equation}
where $M(x)$ is slowly varying as $x\rightarrow 0^+$. Let $g$ be the inverse of $U(1-\cdot)$, and define
\[\Theta(t):=-\int_{0}^t\log f'(1-g(s))ds,\hspace{20pt}\forall t\geq 1.\]
Noting that the size-biasing result at (\ref{sizebias}) yields $\mathbf{E}(s^{Z_n^*})=sf'_n(s)=s\prod_{m=0}^{n-1}f'(f_m(s))$, we are able to proceed as in the proof of \cite{Pakes}, Theorem 2, to deduce that
\[\mathbf{E}\left(s^{Z_n^*}\right)=s\Delta_n(s)\exp\left\{\Theta(U(s))-\Theta(n+U(s))\right\},\hspace{20pt}\forall s\in[0,1),\]
where $\Delta_n(s)\leq 1$. Furthermore, as computed in \cite{Pakes}, the asymptotic behaviour of $U$ described at (\ref{utail}) implies that $\Theta(t)=\alpha(\alpha-1)^{-1}\log t +r(t)$, where the remainder term satisfies $r'(t)=o(t^{-1})$ as $t\rightarrow \infty$. In particular, it follows that
\begin{equation}\label{tuy}
\mathbf{E}\left(s^{Z_n^*}\right)\leq\left(1+\frac{n}{U(s)}\right)^{-\tfrac{\alpha}{\alpha-1}}\exp\left\{r(U(s))-r(n+U(s))\right\},\hspace{20pt}\forall s\in[0,1).
\end{equation}
By definition, $U$ is an increasing function and therefore if $\lambda\in[0,1]$, then $U(e^{-\lambda p_n})\geq U(e^{-p_n})\geq U(1-p_n)=n$, where the final equality is easily checked by applying (\ref{usum}) and the fact that $p_n=1-f_n(0)$. Consequently, given $\eta\in(0,1)$, since $r'(t)=o(t^{-1})$, there exists an integer $n_0$ such that
\[r(U(e^{-\lambda p_n}))-r(n+U(e^{-\lambda p_n}))\leq {\eta a_n(\lambda)},\hspace{20pt}\forall n\geq n_0,\lambda\in[0,1],\]
where we define $a_n(\lambda):=n/U(e^{-\lambda p_n})$. Letting $c_2$ be a constant such that $e^x\leq 1+c_2 x$ for $x\in[0,1]$, then the above inequality and the bound at (\ref{tuy}) imply that
\[\mathbf{E}\left(e^{-\lambda p_nZ_n^*}\right)\leq \frac{1+c_2\eta a_n(\lambda)}{\left(1+a_n(\lambda)\right)^{\tfrac{\alpha}{\alpha-1}}},\hspace{20pt}\forall n\geq n_0,\lambda\in[0,1].\]
Since $c_2$ is independent of the choice of $\eta$, if we are given $\varepsilon'\in(0,\tfrac{\alpha}{\alpha-1})$, then for small enough $\eta$, we have that $1+c_2\eta x\leq (1+x)^{\varepsilon'}$ for every $x\in[0,1]$, and therefore
\[\mathbf{E}\left(e^{-\lambda p_nZ_n^*}\right)\leq \left(1+a_n(\lambda)\right)^{-\left(\tfrac{\alpha}{\alpha-1}-\varepsilon'\right)}\leq 1-c_3a_n(\lambda),\hspace{20pt}\forall n\geq n_0,\lambda\in[0,1],\]
for some constant $c_3$. Thus, to complete the proof of (\ref{gen}), it remains to obtain a suitable lower bound for $a_n(\lambda)$. It follows from (\ref{utail}) and the monotonicity of $U$ that
\[a_n(\lambda)=\frac{n}{U(e^{-\lambda p_n})}\geq \frac{U(1-p_n)}{U(1-c_4\lambda p_n)}=c_5\lambda^{\alpha-1}\frac{M(c_4\lambda p_n)}{M(p_n)}\geq c_6\lambda^{\alpha-1+\varepsilon},\hspace{10pt}\forall n\geq n_0,\lambda\in[0,1],\]
for suitably chosen constants $c_4,c_5,c_6$, where to deduce the final inequality we use the representation theorem for slowly varying functions (see \cite{sen}, Theorem 1.2, for example). This completes the proof of (\ref{gen}).

For any non-negative random variable $\xi$ we have that
\[1-\mathbf{E}\left(e^{-\theta\xi}\right)=\int_0^\infty\mathbf{P}(\xi\geq x)\theta e^{-\theta x}dx,\hspace{20pt}\forall \theta> 0,\]
from which it easily follows that
\[1-\mathbf{E}\left(e^{-\theta\xi}\right)\leq x+\mathbf{P}(\xi\geq x/\theta),\hspace{20pt}\forall \theta,x>0.\]
For $\beta\in(0,1)$, taking $\xi=p_nZ^*_n$, $\theta=\lambda^{-1/(1-\beta)}$ and $x=\lambda \theta$ in the above inequality, we obtain from (\ref{gen}) that
\[\mathbf{P}\left(Z_n^*\geq \lambda p_n^{-1}\right)\geq c_1\lambda^{-(\alpha-1+\varepsilon)/(1-\beta)} - \lambda^{-\beta/(1-\beta)},\hspace{20pt}\forall n\geq n_0,\lambda\geq 1.\]
Now, assume that $\alpha\in(1,2)$ and $\beta_2>(\alpha-1)/(2-\alpha)$. By setting $\beta=\alpha-1+2\varepsilon$ for $\varepsilon$ chosen suitably small, the result follows.
\end{proof}

To prove a polynomial lower bound for the tail of the distribution of $Y^*_n$ near infinity, we will use the fact that $(p_n)_{n\geq 0}$ is regularly varying as $n\rightarrow\infty$, which follows from (\ref{probdecay}).

{\lem \label{pnprop} We can write $p_n=n^{-\frac{1}{\alpha-1}}\ell(n)$, where $\ell(n)$ is a slowly varying function as $n\rightarrow \infty$. Moreover, if $\varepsilon>0$, then there exist constants $c_1,c_2\in(0,\infty)$ such that
\[c_1\left(\frac{n}{m}\right)^{-\varepsilon}\leq \frac{\ell (n)}{\ell(m)}\leq c_2\left(\frac{n}{m}\right)^{\varepsilon},\]
whenever $1\leq m\leq n$.}
\begin{proof} That $p_n=n^{-\frac{1}{\alpha-1}}\ell(n)$, where $\ell(n)$ is a slowly varying function as $n\rightarrow \infty$, follows from (\ref{probdecay}) by applying $5^o$ of \cite{sen}, Section 1.5. The remaining claim can be proved using the representation theorem for slowly varying functions (see \cite{sen}, Theorem 1.2, for example).
\end{proof}

We will also apply the subsequent adaptation of \cite{BarKum}, Lemma 2.3(a), which establishes the result in the case when the offspring distribution is binomial.

{\lem \label{binlem} There exist strictly positive constants $c_1$ and $c_2$ such that
\[\mathbf{P}\left(Y_{2n}\geq c_1 n p_n^{-1}\right)\geq c_2 p_n,\hspace{20pt}\forall n\in \mathbb{N}.\]}
\begin{proof} First observe that
\begin{eqnarray*}
1+2n=\mathbf{E}(Y_{2n})&=&\mathbf{E}(Y_{2n}\mathbf{1}_{\{Z_n=0\}})+\mathbf{E}(Y_{2n}\mathbf{1}_{\{Z_n>0\}})\\
&\leq & \mathbf{E}(Y_{n})+p_n\mathbf{E}(Y_{2n}|Z_n>0)\\
&=&n+1 +p_n\mathbf{E}(Y_{2n}|Z_n>0).
\end{eqnarray*}
In particular, this implies that
\[np_n^{-1}\leq \mathbf{E}(Y_{2n}|Z_n>0).\]
Furthermore, if $\beta\in (0,\alpha-1)$, then
\[\mathbf{E}(Y_{2n}^{1+\beta}|Z_n>0)\leq p_n^{-1}\mathbf{E}(Y_{2n}^{1+\beta})
\leq p_n^{-1}(2n+1)^{1+\beta}\mathbf{E}\left(\max_{m\leq 2n}Z_m^{1+\beta}\right).\]
Since $(Z_n)_{n\geq 0}$ is a martingale, we can apply Doob's martingale norm inequality to obtain from this that
\[\mathbf{E}(Y_{2n}^{1+\beta}|Z_n>0)\leq \left(\frac{1+\beta}{\beta}\right)^{1+\beta} p_n^{-1}(2n+1)^{1+\beta}\mathbf{E}\left(Z_{2n}^{1+\beta}\right)\leq c_1\left(np_{n}^{-1}\right)^{1+\beta},\]
for some finite constant $c_1$, where we apply Lemma \ref{moments} to bound $\mathbf{E}(Z_{2n}^{1+\beta})$ by $c_2p_{2n}^{-\beta}$ and Lemma \ref{pnprop} to show that $p_{2n}^{-\beta}\leq c_3p_n^{-\beta}$.

Now, let $\varepsilon\in (0,1)$ and $\xi$ be a non-negative random variable, then by H\"{o}lder's inequality we have that
\begin{equation}\label{revhold}
(1-\varepsilon)\mathbf{E}(\xi)\leq \mathbf{E}\left(\xi\mathbf{1}_{\{\xi\geq \varepsilon\mathbf{E}(\xi)\}}\right)\leq \mathbf{E}\left(\xi^{1+\beta}\right)^{1/(1+\beta)}\mathbf{P}\left(\xi\geq \varepsilon\mathbf{E}(\xi)\right)^{\beta/(1+\beta)},
\end{equation}
assuming that the appropriate moments are finite. Applying this bound to $Y_{2n}$ with respect to the conditioned measure $\mathbf{P}(\cdot|Z_n>0)$, the above estimates yield
\[\mathbf{P}\left(Y_{2n}\geq \varepsilon np_n^{-1}|Z_n>0\right)\geq c_4>0,\hspace{20pt}\forall n\in\mathbb{N},\]
for some constant $c_4$. Hence, we have that
\[\mathbf{P}\left(Y_{2n}\geq \varepsilon n p_n^{-1}\right)\geq p_n\mathbf{P}\left(Y_{2n}\geq \varepsilon n p_n^{-1}|Z_n>0\right)\geq c_4p_n,\hspace{20pt}\forall n\in\mathbb{N},\]
which completes the proof.
\end{proof}

We can now prove our first tail bounds for $Y_n^*$. To obtain the upper polynomial tail bound near infinity, we apply the size-biased interpretation of the law of the first $n$ generations of $T^*$ and a standard martingale bound. In the proof of the corresponding lower bound, we rely on decomposition of $T^*$ that appears in \cite{Kesten}. The same decomposition will also be applied in Proposition \ref{ylower} and Lemma \ref{upc} below. Henceforth, we will use the notation $\mathrm{Bin}(N,p)$ to represent a binomial random variable with parameters $N$ and $p$.

{\propn\label{yupper}  For $\beta_1\in (0,\alpha-1)$, there exists a finite constant $c_1$ such that
\begin{equation}\label{yupperb}
\mathbf{P}\left(Y_n^*\geq \lambda np_n^{-1}\right)\leq c_1 \lambda^{-\beta_1},\hspace{20pt}\forall n\in \mathbb{N},\lambda>0.
\end{equation}
Moreover, for $\alpha\in(1,2)$, $\beta_2>(\alpha-1)/(2-\alpha)$, there exists a strictly positive constant $c_2$ and integer $n_0$ such that
\begin{equation}\label{yupperc}
\mathbf{P}\left(Y_n^*\geq \lambda np_n^{-1}\right)\geq c_2 \lambda^{-\beta_2},\hspace{20pt}\forall n\geq n_0,\lambda\geq 1.
\end{equation}}
\begin{proof} Fix $\beta_1\in (0,\alpha-1)$. It follows from the size-biasing result at (\ref{sizebias}) that
\[\mathbf{E}\left({Y_n^*}^{\beta_1}\right)= \mathbf{E}\left(Y_n^{\beta_1} Z_n\right)\leq (n+1)^{\beta_1}\mathbf{E}\left(\max_{m\leq n}Z_m^{1+\beta_1}\right).\]
Applying Doob's martingale norm inequality and Lemma \ref{moments}, we consequently obtain that there exists a finite constant $c_1$ such that
\begin{equation}\label{ymoments}
\mathbf{E}\left({Y_n^*}^{\beta_1}\right)\leq c_1 (np_n^{-1})^{\beta_1}.
\end{equation}
The result at (\ref{yupperb}) is readily deduced from this bound.

Now assume that $\alpha\in(1,2)$, $\beta_2>(\alpha-1)/(2-\alpha)$, and let $c_1$ and $c_2$ be the constants of Lemma \ref{binlem}. Clearly, we have that
\[\mathbf{P}\left(Y^*_{3n}\geq \lambda np_n^{-1}\right)\geq \mathbf{P}\left(Y^*_{3n}\geq \lambda np_n^{-1}|Z^*_{n}\geq c_3 \lambda p_n^{-1}\right)\mathbf{P}\left(Z^*_{n}\geq c_3\lambda p_n^{-1}\right),\hspace{10pt}\forall n\in\mathbb{N},\lambda\geq 1,\]
for an arbitrary constant $c_3\geq 2$. By \cite{Kesten}, Lemma 2.2, the tree $T^*$ has a unique infinite line of descent, or backbone, and the descendants of the individuals in the $n$th generation of $T^*$ which are not on the backbone have the same distribution as the unconditioned $T$, independently of each other; hence the first factor above is bounded below by $\mathbf{P}(Y_{2n}[\lfloor c_3 \lambda p_n^{-1}\rfloor-1] \geq \lambda np_n^{-1})$, where $Y_{2n}[m]$ is the sum of $m$ independent copies of $Y_{2n}$. Thus, applying Lemma \ref{binlem}, we obtain
\[\mathbf{P}\left(Y^*_{3n}\geq \lambda np_n^{-1}|Z^*_{n}\geq c_3 \lambda p_n^{-1}\right)\geq  \mathbf{P}\left(\mathrm{Bin}\left(\lfloor c_3 \lambda p_n^{-1}\rfloor-1, c_2p_n\right)\geq c_1^{-1} \lambda\right).\]
Taking $c_3$ large enough, the ``reverse H\"{o}lder'' inequality of (\ref{revhold}) implies that the right-hand side is bounded below by a strictly positive constant $c_4$ uniformly in $n\in\mathbb{N}$ and $\lambda\geq 1$. Consequently, by (\ref{zupperc}),
\[\mathbf{P}\left(Y^*_{3n}\geq \lambda np_n^{-1}\right)\geq c_4\mathbf{P}\left(Z^*_{n}\geq c_3\lambda p_n^{-1}\right)\geq c_5\lambda^{-\beta_2},\hspace{20pt}\forall n\geq n_0,\lambda\geq 1.\]
From this we can deduce (\ref{yupperc}) by applying the monotonicity of $(Y^*_n)_{n\geq 0}$ and Lemma \ref{pnprop}.
\end{proof}

We now consider the tail near 0 of the distributions of the random variables $Z^*_n$. In particular, to deduce a polynomial upper bound, we follow a generating function argument in which we apply the known asymptotics of the sequence of survival probabilities $(p_n)_{n\geq 0}$.

{\propn \label{zlower} For $\beta\in (0,\alpha-1)$, there exists a finite constant $c$ such that
\[\mathbf{P}\left(Z_n^*\leq \lambda p_n^{-1}\right)\leq c \lambda^{\beta},\hspace{20pt}\forall n\in \mathbb{N},\lambda>0.\]}
\begin{proof} Fix $\beta\in(0,\alpha-1)$. We will start by showing that there exists a finite constant $c_1$ such that
\begin{equation}\label{cf}
\mathbf{E}\left(e^{-\lambda p_nZ_n}| Z_n>0\right)\leq c_1\lambda^{-\beta},\hspace{20pt}\forall n\in \mathbb{N},\lambda\in[1,p_n^{-1}].
\end{equation}
Clearly, we have that
\[\mathbf{E}\left(e^{-\lambda p_nZ_n}| Z_n>0\right)=1-\frac{1-\mathbf{E}\left(e^{-\lambda p_nZ_n}\right)}{p_n}.\]
Choose an integer $k=k(n,\lambda)\geq 0$ as in the proof of \cite{Slack}, Theorem 1, to satisfy $p_k\geq 1-e^{-\lambda p_n}>p_{k+1}$, then, by the Markov property of $(Z_n)_{n\geq 0}$,
\[\mathbf{E}\left(e^{-\lambda p_nZ_n}\right)\leq \mathbf{E}\left((1-p_{k+1})^{Z_n}\right)=1-p_{n+k+1}.\]
Hence,
\[\mathbf{E}\left(e^{-\lambda p_nZ_n}| Z_n>0\right)\leq 1-\frac{p_{n+k+1}}{p_n}.\]
In the proof of \cite{Slack}, Theorem 1, it is observed that
\[p_{m+1}p_m^{-1}\geq 1- c_2m^{-1},\hspace{20pt}\forall m\in\mathbb{N},\]
for some finite constant $c_2$. Applying this bound and the inequality $(1-x)^n\geq 1-nx$ for every $x\in[0,1]$ and $n\in \mathbb{N}$, it is possible to deduce the existence of a finite constant $c_3$ such that
\[\mathbf{E}\left(e^{-\lambda p_nZ_n} | Z_n>0\right)\leq c_3 (k+1)n^{-1},\]
for every $n\in \mathbb{N}$ and $\lambda\geq 1$. To estimate $(k+1)n^{-1}$, we first choose $c_4$ small enough so that $e^{-x}\leq 1-c_4x$ for $x\in[0,1]$, which implies $p_k\geq c_4 \lambda p_n$ for every $n\in \mathbb{N}$ and $\lambda\in [1,p_n^{-1}]$. This inequality allows us to apply Lemma \ref{pnprop} to demonstrate that there exists a finite constant $c_5$ such that $k+1\leq c_5\lambda^{-\beta}n$ for every $n\in \mathbb{N}$ and $\lambda\in [1,p_n^{-1}]$, which completes the proof of (\ref{cf}).

Before continuing, note that a simple coupling argument allows us to obtain that
\[\mathbf{P}\left(Z_{n+m}>0|Z_n\in(0, \lambda]\right)\leq \mathbf{P}\left(Z_{n+m}>0|Z_n>0\right),\]
for any $m,n\in\mathbb{N}$ and $\lambda>0$. By Bayes' formula, this is equivalent to
\[\mathbf{P}\left(Z_n\leq \lambda | Z_{m+n}>0\right)\leq \mathbf{P}\left(Z_n\leq \lambda | Z_n>0\right),\]
for any $m,n\in \mathbb{N}$ and $\lambda>0$. Thus,
\begin{eqnarray*}
\mathbf{P}\left(Z^*_n\leq \lambda p_n^{-1}\right)&=&\lim_{m\rightarrow\infty}\mathbf{P}\left(Z_n\leq \lambda p_n^{-1} | Z_{m+n}>0\right)\\
&\leq & \mathbf{P}\left(Z_n\leq \lambda p_n^{-1} | Z_{n}>0\right)\\
&\leq & \mathbf{E}\left(e^{1-\lambda^{-1} p_nZ_n} | Z_n>0\right)\\
&\leq & c_6\lambda^{\beta},
\end{eqnarray*}
whenever $n\in\mathbb{N}$ and $\lambda\in[p_n,1]$. Since the claim of the proposition is trivial for $\lambda<p_n$ and $\lambda>1$, the proof is complete.
\end{proof}

This result allows us to prove a tail bound near 0 for $Y^*_n$ that is uniform in $n$.

{\propn\label{ylower} For $\gamma\in (0,1-\alpha^{-1})$, there exists a finite constant $c$ such that
\[\mathbf{P}\left(Y^*_n\leq \lambda n p_n^{-1}\right)\leq c \lambda^{\gamma},\hspace{20pt}\forall n\in \mathbb{N}, \lambda>0.\]}
\begin{proof} Fix $\gamma\in (0,1-\alpha^{-1})$ and choose $\beta\in (0,\alpha-1)$ large enough so that $\gamma'=\gamma/\beta<\alpha^{-1}$. Let $c_1$ and $c_2$ be the constants of Lemma \ref{binlem}. We will prove the result for $\lambda\in [p_n, c_1]$, from which the result for any $\lambda>0$ follows easily. We can write
\[\mathbf{P}\left(Y^*_{3n}\leq \lambda np_n^{-1}\right)\leq \mathbf{P}\left(Z^*_{n}\leq \lambda^{\gamma'} p_n^{-1}\right) +\mathbf{P}\left(Y^*_{3n}\leq \lambda np_n^{-1},Z^*_{n}> \lambda^{\gamma'} p_n^{-1}\right).\]
By Proposition \ref{zlower}, there exists a finite constant $c_3$ such that the first term here is bounded above by $c_3\lambda^{\gamma}$ for any $n\in\mathbb{N}$ and $\lambda>0$. By applying the decomposition of $T^*$ described in the proof of Proposition \ref{yupper}, we have that the second term is bounded above by $\mathbf{P}(Y_{2n}[\lfloor \lambda^{\gamma'}p_n^{-1}\rfloor]\leq \lambda np_n^{-1})$, where $Y_{2n}[m]$ is the sum of $m$ independent copies of $Y_{2n}$. If we choose $m=m(n,\lambda)$ to be the smallest integer such that $\lambda n p_n^{-1}<c_1mp_{m}^{-1}$, then $m\leq n$ and, applying Lemma \ref{binlem}, we obtain that
\begin{eqnarray*}
\mathbf{P}\left(Y^*_{3n}\leq \lambda np_n^{-1},Z^*_{n}> \lambda^{\gamma'} p_n^{-1}\right)&\leq&\mathbf{P}\left(c_1mp_m^{-1}\mathrm{Bin}(\lfloor \lambda^{\gamma'}p_n^{-1}\rfloor, c_2p_m)\leq \lambda np_n^{-1}\right)\\
&\leq &\mathbf{P}\left(\mathrm{Bin}(\lfloor \lambda^{\gamma'}p_n^{-1}\rfloor, c_2p_m)<1\right)\\
&=&(1-c_2p_m)^{\lfloor \lambda^{\gamma'}p_n^{-1}\rfloor}\\
&\leq & e^{-c_2p_m \lfloor \lambda^{\gamma'}p_n^{-1}\rfloor}.
\end{eqnarray*}
It is an elementary exercise to apply Lemma \ref{pnprop} to deduce that our choice of $m$ implies that if $\gamma''\in(\gamma',\alpha^{-1})$, then there exists a constant $c_4>0$ such that
$p_mp_n^{-1}\geq c_4 \lambda^{-\gamma''}$ for every $n\in\mathbb{N}$ and $\lambda \in [p_n,c_1]$. Consequently,
\[\mathbf{P}\left(Y^*_{3n}\leq \lambda np_n^{-1}\right)\leq c_3\lambda^{\gamma} + e^{-c_5 \lambda^{(\gamma'-\gamma'')}}\leq c_6\lambda^{\gamma},\hspace{20pt}\forall n\in\mathbb{N}, \lambda\in [p_n,c_1],\]
from which the result follows by applying the monotonicity of $(Y^*_n)_{n\geq 0}$ and Lemma \ref{pnprop}.
\end{proof}

Finally, we prove a tail bound for the number of individuals in the $n$th generation of $T^*$ that have descendants in the $2n$th generation, which we denote by $M_n^{2n}$.

{\lem \label{upc} For every $\beta\in(0,\alpha-1)$, there exists a finite constant $c$ such that
\[\mathbf{P}\left(M_n^{2n}\geq \lambda\right)\leq c\lambda^{-\beta},\hspace{20pt}\forall n\in\mathbb{N}, \lambda>0.\]}
\begin{proof} Fix $\beta\in(0,\alpha-1)$. If we condition on the first $n$ generations of $T^*$, denoted by $T^*|_n$, and the backbone $B$, then \cite{Kesten}, Lemma 2.2 implies that
\[\mathbf{P}\left(M_n^{2n}\geq \lambda~|~\,  T^*|_n, B\right)
%\[\mathbf{P}\left(M_n^{2n}\geq \lambda|  T^*|_n, B\right)
=\mathbf{P}\left(\mathrm{Bin}(N,p_n)\geq \lambda-1\right)|_{N=Z^*_n-1}.\]
Consequently,
\[\mathbf{P}\left(M_n^{2n}\geq \lambda\right)\leq \mathbf{P}\left(p_nZ^*_n\geq \lambda/2\right) +\mathbf{P}\left(\mathrm{Bin}(\lceil\tfrac{\lambda}{2p_n}\rceil,p_n)\geq \lambda-1\right).\]
Thus, Proposition \ref{zupper} and Chebyshev's inequality imply that there exists a finite constant $c$ such that
\[\mathbf{P}\left(M_n^{2n}\geq \lambda\right)\leq c\lambda^{-\beta} +\frac{\lceil\tfrac{\lambda}{2p_n}\rceil p_n}{(\lambda-1-\lceil\tfrac{\lambda}{2p_n}\rceil p_n)^2}.\] The result follows.
\end{proof}

\section{Proof of initial random walk results}\label{proofs}

In this section we complete the proofs of Theorems \ref{probthm}, \ref{annthm}, \ref{qthm} and \ref{exponents}, though we first introduce some further notation that we will apply. The volume of a ball of radius $R$ about the root of $T^*$ is given by
\[V(R):=\mu^{T^*}(B(R)),\]
where $B(R):=\{x\in T^*:d_{T^*}(\rho,x)\leq R\}$ and $\mu^{T^*}$ is the invariant measure of $X$ defined in the introduction. Let $\mathcal{E}$ be a quadratic form on $\mathbb{R}^{T^*}$ that satisfies
\[\mathcal{E}(f,g)=\frac{1}{2}\sum_{\substack{x,y\in T^*\\x\sim y}}(f(x)-f(y))(g(x)-g(y)),\]
where $x\sim y$ if and only if $\{x,y\}$ is an edge of $T^*$. The quantity $\mathcal{E}(f,f)$ represents the energy dissipation when we suppose that $T^*$ is an electrical network with a unit resistor placed along each edge and vertices are held at the potential $f$. The associated effective resistance operator can be defined by
\[R_{eff}(A,B)^{-1}:=\inf\{\mathcal{E}(f,f):f|_A=1,f|_B=0\},\]
for disjoint subsets $A,B\subseteq T^*$.

Recall the volume growth function $v$ defined at (\ref{vgf}), and let $r:\mathbb{R}_+\rightarrow\mathbb{R}_+$ be the identity function on $\mathbb{R}_+$. It is clear that $v$ and $r$ are both strictly increasing functions, $r$ satisfies $r(R)/r(R')=R/R'$ for every $1\leq R'\leq R<\infty$, and, by applying Lemma \ref{pnprop}, we can check that for each $\varepsilon>0$ there exist constants $c_1,c_2\in(0,\infty)$ such that
\[c_1\left(\frac{R}{R'}\right)^{\frac{\alpha}{\alpha-1}-\varepsilon}\leq\frac{v(R)}{v(R')}\leq c_2\left(\frac{R}{R'}\right)^{\frac{\alpha}{\alpha-1}+\varepsilon},\]
whenever $1\leq R'\leq R<\infty$. Consequently, the conditions required on $v$ and $r$ in \cite{KM} are fulfilled (in \cite{KM}, it was also assumed that $v(1)=r(1)=1$, but we can easily neglect this technicality by adjusting constants as necessary), and to deduce many of the results about the random walk on $T^*$ stated in the introduction, it will suffice to check that the relevant parts of \cite{KM}, Assumption 1.2 are satisfied. More specifically, we will check that, if we denote by
\[J(\lambda):=\{R\in[1,\infty]:\lambda^{-1}v(R)\leq V(R)\leq \lambda v(R), R_{eff}(\{\rho\},B(R)^c)\geq\lambda^{-1}r(R)\},\]
for $\lambda\geq 1$, then the probability that $R\in J(\lambda)$ is bounded below, uniformly in $R$, by a function of $\lambda$ that increases to 1 polynomially. This result explains why $v$ can be thought of as a volume growth function for $T^*$. Note that, in \cite{KM}, $J(\lambda)$ has the extra restriction that $R_{eff}(\rho,x)\leq \lambda r(d_{T^*}(\rho,x))$ for every $x\in B(R)$. However, since $T^*$ is a tree, this condition always holds, and so we omit it.

{\lem \label{ass12} $T^*$ satisfies Assumptions 1.2(1) and 1.2(3) of \cite{KM}. In particular, for every $\gamma\in(0,1-\alpha^{-1})$, there exists a finite constant $c$ such that
\[\inf_{R\geq 1}\mathbf{P}\left(R\in J(\lambda)\right)\geq 1-c\lambda^{-\gamma},\hspace{20pt}\forall \lambda\geq1.\]}
\begin{proof} Fix $\gamma\in(0,1-\alpha^{-1})$. First note that, since $T^*$ is a tree, we have that
\begin{equation}\label{215}
Y_R^*\leq V(R)\leq 2Y_{R+1}^*,\hspace{20pt}\forall R\in \mathbb{N}.
\end{equation}
Therefore it will be adequate to prove the result for $R\in\mathbb{N}$ and $V(R)$ replaced by $Y_R^*$. That
\[\inf_{R\in\mathbb{N}}\mathbf{P}\left(\lambda^{-1}v(R)\leq Y_R^*\leq \lambda v(R)\right)\geq1-c_1\lambda^{-\gamma},\hspace{20pt}\forall \lambda\geq 1,\]
for some finite constant $c_1$ is an easy consequence of Propositions \ref{yupper} and \ref{ylower}. By imitating the proof of \cite{BarKum}, Lemma 4.5, it is possible to prove that
\[R_{eff}(\rho,B(2R)^c)\geq \frac{R}{M_{R}^{2R}},\]
for every $R\in \mathbb{N}$. Thus, applying Lemma \ref{upc}, we have that
\[\inf_{R\geq 1}\mathbf{P}\left(R_{eff}(\rho,B(R)^c)\geq\lambda^{-1}r(R)\right)\geq1-c_2\lambda^{-\gamma},\hspace{20pt}\forall \lambda\geq 1,\]
for some finite constant $c_2$, and the lemma holds as claimed.
\end{proof}

\begin{proof}[Proof of Theorem \ref{probthm}] Apart from (\ref{kestens}), the limits can all be obtained using \cite{KM}, Proposition 1.3. Since $d_{T^*}(\rho, X_m)\geq R$ implies that $\tau_R\leq m$, the right-hand inequality of (\ref{kestens}) follows from the left-hand inequality of (\ref{dlim}). Consequently it remains to show that
\[\lim_{\lambda\rightarrow\infty}\inf_{R\in\mathbb{N}}\mathbb{P}\left(\lambda^{-1}h(R)\leq \tau_R\right)=1.\]
By \cite{KM}, Proposition 3.5(a), there exist constants $c_1,c_2,c_3,c_4\in(0,\infty)$ that depend only $\lambda$ such that, if $\varepsilon<c_1$ and $R,\varepsilon R, c_2\varepsilon R\in J(\lambda)$, then
\[P^{T^*}_\rho\left(\tau_R\leq c_3\varepsilon^\beta h(R)\right)\leq c_4 \varepsilon,\]
for some deterministic constant $\beta>0$. Hence, for any $\lambda>0$,
\begin{eqnarray*}
\lim_{\varepsilon\rightarrow 0}\limsup_{R\rightarrow\infty}\mathbb{P}\left(\tau_R\leq \varepsilon h(R)\right)&=&\lim_{\varepsilon\rightarrow 0}\limsup_{R\rightarrow\infty}\mathbb{P}\left(\tau_R\leq c_3\varepsilon^\beta h(R)\right)\\
&\leq&\lim_{\varepsilon\rightarrow 0}\limsup_{R\rightarrow\infty}\left\{c_4\varepsilon+1 - \mathbf{P}\left(R,\varepsilon R, c_2\varepsilon R\in J(\lambda)\right)\right\}\\
&\leq & 3 \sup_{R\geq 1}\mathbf{P}\left(R\not\in J(\lambda)\right).
\end{eqnarray*}
Since $\lambda$ is arbitrary, we can make this upper bound as small as we choose by applying Lemma \ref{ass12}, and so $\lim_{\varepsilon\rightarrow 0}\limsup_{R\rightarrow\infty}\mathbb{P}(\tau_R\leq \varepsilon h(R))=0$. The desired conclusion is readily deduced from this limit.
\end{proof}

\begin{proof}[Proof of Theorem \ref{annthm}] The proofs of the lower bounds at (\ref{anntau}), (\ref{anntd}) and (\ref{annd}) require only straightforward adaptations of the proofs of the lower bounds in \cite{KM}, Proposition 1.4, and are omitted. As in \cite{BarKum}, Proposition 4.4, for example, we have that $E^{T^*}_\rho \tau_R\leq (R+1)V(R)$. Since $V(R)\leq 2Y_{R+1}^*$ (see (\ref{215})), it follows that
\[\mathbf{E}\left(\left(E^{T^*}_\rho\tau_R\right)^\beta\right)\leq 2^\beta(R+1)^\beta\mathbf{E}\left({Y^*_{R+1}}^\beta\right),\hspace{20pt}\forall R\in\mathbb{N}.\]
Thus the upper bound at (\ref{anntau}) follows from the estimate of the $\beta$-moments of $(Y^*_n)_{n\geq 0}$ that appears at (\ref{ymoments}).

Similarly to the proof of \cite{KM}, Remark 1.6.1, it is possible to deduce that there exists a finite constant $c_1$ such that
\begin{equation}\label{hkupper}
\mathbf{E}\left(p_{2m}^{T^*}(\rho,\rho)^\gamma\right)\leq \frac{c_1}{v(\mathcal{I}(m))^\gamma}\mathbf{E}\left(1+\frac{v(R)^\gamma}{{Y^*_R}^\gamma}\right),\hspace{20pt}\forall m\in\mathbb{N},
\end{equation}
where $R=R(m)$ is chosen to satisfy $\tfrac{1}{2}h(R)\leq m\leq h(R)$, and we have again applied (\ref{215}). By the tail bound of Proposition \ref{ylower}, if $\gamma$ is in the range $(0,1-\alpha^{-1})$, then we can bound the expectation on the right-hand side of the above expression uniformly in $R$ by a constant. This completes the proof of (\ref{anntd}).

It remains to prove the upper bound at (\ref{annd}). First, let $\tilde{\tau}_R$ be the first hitting time of the vertex on the backbone at a distance $R$ from the origin, i.e.
\begin{equation}\label{tildetau}
\tilde{\tau}_R:=\min\{n:X_n\in B,d_{T^*}(\rho,X_n)=R\},
\end{equation}
where $B\subseteq T^*$ is the backbone of $T^*$ (the unique non-intersecting infinite path in $T^*$ which starts at the root $\rho$). By (\ref{kestens}), it is clear that $\varepsilon>0$ can be chosen small enough so that $\mathbb{P}\left(\tilde{\tau}_R\leq \varepsilon h(R)\right)\leq \mathbb{P}\left({\tau}_R\leq \varepsilon h(R)\right)\leq \tfrac{1}{2}$ for every $R\in\mathbb{N}$, which implies that
\[\mathbb{P}\left(\tilde{\tau}_R\leq t\right)\leq \frac{1}{2}+\frac{t}{\varepsilon h(R)},\hspace{20pt}\forall R\in \mathbb{N}, t>0.\]
Now observe that if we define, for $i,R\in\mathbb{N}$,
\[\tilde{\tau}_{R}^{i}:=\#\{n\in[\tilde{\tau}_{(i-1)R},\tilde{\tau}_{iR}):X_n,X_{n+1}\mbox{ equal $b_{(i-1)R}$ or are descendents of }b_{(i-1)R}\},\]
where $b_{(i-1)R}$ is the vertex on the backbone at a distance $(i-1)R$ from $\rho$, then with respect to the annealed measure $\mathbb{P}$ the elements of the sequence $(\tilde{\tau}_{R}^{i})_{i\geq 1}$ are independent and have the same distribution as $\tilde{\tau}_R$ (this follows from the description of $T^*$ given in Lemma 2.2 of \cite{Kesten}). Thus, since $\tilde{\tau}_{nR}\geq \sum_{i=1}^n \tilde{\tau}_{R}^{i}$, we can apply Lemma 1.1 of \cite{BarBas} to obtain that
\[\log\mathbb{P}\left(\tilde{\tau}_{nR}\leq t\right)\leq 2\left(\frac{2nt}{\varepsilon h(R)}\right)^{1/2}-n\log 2,\hspace{20pt}\forall n,R\in\mathbb{N}, t>0.\]
In particular, by setting $t=c_2nh(R)$ for constant $c_2$ chosen suitably small, it follows that
\begin{equation}\label{tail}
\mathbb{P}\left(\tilde{\tau}_{nR}\leq c_2nh(R)\right)\leq e^{-c_3n},\hspace{20pt}\forall n,R\in\mathbb{N},
\end{equation}
where $c_3$ is a strictly positive constant.

For $m\in \mathbb{N}$, write $R=\lfloor\mathcal{I}(m)\rfloor$, then, for every $\lambda\in\mathbb{N}$, $\varepsilon\in[0,1]$ and $\eta>0$,
\begin{eqnarray}
\mathbb{P}\left(\max_{0\leq k\leq m}d_{T^*}(\rho,X_k)\geq \lambda \mathcal{I}(m)\right)&\leq&\mathbb{P}\left(\tau_{\lambda R}\leq m\right)\nonumber\\
&\leq & \mathbb{P}\left(\tau_{\lambda R}\leq m,M_{\lfloor\varepsilon\lambda R\rfloor}^{\lambda R}=1\right)\nonumber\\
&&+\mathbf{P}\left(M_{\lfloor\varepsilon\lambda R\rfloor}^{\lambda R}>1, Z^*_{\lfloor\varepsilon\lambda R\rfloor}\leq \eta p_{\lfloor\varepsilon\lambda R\rfloor}^{-1}+1\right)\nonumber\\
&&+\mathbf{P}\left(Z^*_{\lfloor\varepsilon\lambda R\rfloor}> \eta p_{\lfloor\varepsilon\lambda R\rfloor}^{-1}+1\right),\label{decomp}
\end{eqnarray}
where, generalising the notation of the previous section, $M_n^{m+n}$ is the number of individuals in the $n$th generation of $T^*$ that have descendants in the $(m+n)$th generation. On the event $\{M_{\lfloor\varepsilon\lambda R\rfloor}^{\lambda R}=1\}$, of the vertices in generation $\lfloor\varepsilon\lambda R\rfloor$, only the one on the backbone has descendants in generation $\lambda R$; thus if $X$ has reached generation $\lambda R$ no later than time $m$, then $X$ must have already visited the vertex on the backbone at a distance $\lfloor\varepsilon\lambda R\rfloor$ from the root. Hence, if $c_2\varepsilon \lambda \geq 1$, then it is possible to check that
\[\mathbb{P}\left(\tau_{\lambda R}\leq m,M_{\lfloor\varepsilon\lambda R\rfloor}^{\lambda R}=1\right)\leq \mathbb{P}\left(\tilde{\tau}_{\lfloor\varepsilon\lambda R\rfloor}\leq m\right)\leq \mathbb{P}\left(\tilde{\tau}_{\lfloor\varepsilon\lambda R\rfloor}\leq c_2\varepsilon\lambda m\right)\leq e^{-c_4\varepsilon\lambda},\]
for some constant $c_4>0$, where we apply the bound at (\ref{tail}) to deduce the final inequality. For the second term at (\ref{decomp}), we proceed similarly to the proof of Lemma \ref{upc} to obtain that
\begin{eqnarray*}
\mathbf{P}\left(M_{\lfloor\varepsilon\lambda R\rfloor}^{\lambda R}>1, Z^*_{\lfloor\varepsilon\lambda R\rfloor}\leq \eta p_{\lfloor\varepsilon\lambda R\rfloor}^{-1}+1\right)&\leq& \mathbf{P}\left(\mathrm{Bin}(\lceil\eta p_{\lfloor\varepsilon\lambda R\rfloor}^{-1}\rceil,p_{\lambda R -  \lfloor\varepsilon\lambda R\rfloor})>0\right)\\
&=&1-\left(1-p_{\lambda R -  \lfloor\varepsilon\lambda R\rfloor}\right)^{\lceil\eta p_{\lfloor\varepsilon\lambda R\rfloor}^{-1}\rceil}\\
&\leq & p_{\lambda R -  \lfloor\varepsilon\lambda R\rfloor} \lceil\eta p_{\lfloor\varepsilon\lambda R\rfloor}^{-1}\rceil.
\end{eqnarray*}
For $\beta\in(0,\alpha-1)$, we can bound the third term of (\ref{decomp}) by $c_5\eta^{-\beta}$ by Proposition \ref{zupper}. Combining these bounds, we have that
\[\mathbb{P}\left(\max_{0\leq k\leq m}d_{T^*}(\rho,X_k)\geq \lambda \mathcal{I}(m)\right)\leq e^{-c_4\varepsilon\lambda} + p_{\lambda R -  \lfloor\varepsilon\lambda R\rfloor} \lceil\eta p_{\lfloor\varepsilon\lambda R\rfloor}^{-1}\rceil+c_5\eta^{-\beta},\]
whenever $c_2\varepsilon \lambda \geq 1$.

Finally, fix $\delta\in(0,\alpha^{-1})$ and let $\delta'\in(\delta,\alpha^{-1})$. Choose $\theta_1\in(0,1)$ and $\beta\in(0,\alpha-1)$ large enough so that $\theta_1\beta(1+\beta)^{-1}(\alpha-1)^{-1}\in(\delta',\alpha^{-1})$, and set $\theta_2=\theta_1(1+\beta)^{-1}(\alpha-1)^{-1}$. In the above argument, if we let $\varepsilon=\lambda^{-\theta_1}$ and $\eta=\lambda^{\theta_2}$, then we see that, for every $m\in\mathbb{N}$,
\begin{eqnarray*}
\mathbb{P}\left(\max_{0\leq k\leq m}d_{T^*}(\rho,X_k)\geq \lambda \mathcal{I}(m)\right)&\leq& e^{-c_4\lambda^{1-\theta_1}}+c_6\lambda^{\theta_2-\tfrac{\theta_1}{\alpha-1}}\frac{\ell(\lfloor\lambda^{1-\theta_1} R\rfloor)}{\ell(\lambda R)}+c_5\lambda^{-\theta_2\beta},\\
&\leq & c_7\lambda^{-\delta'},
\end{eqnarray*}
whenever $c_2\lambda^{1-\theta_1}\geq 1$, for some finite constant $c_7$. Note that, to deduce the inequalities above, we have applied the description of the non-extinction probabilities provided by Lemma \ref{pnprop}. Clearly, by increasing $c_7$ if necessary, we can extend this bound to hold for any $\lambda>0$. The upper bound at (\ref{annd}) is an easy consequence of this result.
\end{proof}

{\rem We note that the upper bound for $\mathbb{P}\left(\max_{0\leq k\leq m}d_{T^*}(\rho,X_k)\geq \lambda \mathcal{I}(m)\right)$ obtained in the above proof also implies that (\ref{dlim}) holds when $d_{T^*}(\rho,X_m)$ is replaced by $\max_{0\leq k\leq m}d_{T^*}(\rho,X_k)$.}

\begin{proof}[Proof of Theorem \ref{qthm}] This is an immediate application of Lemma \ref{ass12} and \cite{KM}, Theorem 1.5.
\end{proof}

\begin{proof}[Proof of Theorem \ref{exponents}] Since for a slowly varying function $\ell$, we have that $\log \ell(n)/\log n$ converges to 0 as $n\rightarrow \infty$ (see \cite{sen}, Section 1.5, for example), the claims at (\ref{exptau}), (\ref{exptd}) and (\ref{expd}) follow from Theorem \ref{qthm} and Lemma \ref{pnprop}. Furthermore, that $\log \mu^{T^*}(W_m)/\log m$ converges to $\alpha/(2\alpha-1)$, $P_x^{T^*}$-a.s. for every $x\in T^*$, for $\mathbf{P}$-a.e. realisation of $T^*$ is also proved in \cite{KM}, Theorem 1.5. Since we know from (\ref{215}) that $Y_R^*\leq V(R)\leq 2Y_{R+1}^*$, the proof of the remaining limit involving $\# W_m$ can be obtained by making only minor changes to the proof of the previous result, and so is omitted.
\end{proof}

\section{Annealed off-diagonal transition density}\label{offd}

In addition to the on-diagonal transition density behaviour that we have already established, it is of interest to determine how the transition density of the simple random walk decays away from the diagonal, and in this section we consider the annealed version of this problem. In \cite{BarKum}, for the case of a binomial offspring distribution, bounds for the expectation of the off-diagonal transition density $p_n^{T^*}(x,y)$, conditional on the tree $T^*$ containing the arguments $x$ and $y$, were proved. However, it seems difficult to transfer the arguments used in \cite{BarKum} to the case of a general offspring distribution, even one with finite variance, as doing so would require good uniform control of the laws of $T^*$ conditioned on containing particular vertices. To avoid this issue, we will study the annealed transition density along the backbone of $T^*$. Since this path is $\mathbf{P}$-a.s. present, no conditioning is required in the bounds we prove.
Throughout this section, we will denote the backbone by $\{\rho=b_0,b_1,b_2,\dots\}$, where $b_r$ is the vertex on the backbone satisfying $d_{T^*}(\rho,b_r)=r$. We start by proving an annealed upper transition density bound along the backbone.

{\propn If $\gamma\in(0,1-\alpha^{-1})$, then there exist constants $c_1,c_2\in(0,\infty)$ such that
\[\mathbf{E}\left(p_{2m}^{T^*}(\rho,b_{2r})^\gamma\right)^{\tfrac{1}{\gamma}}\leq c_1 v(\mathcal{I}(m))^{-1}\exp\left\{-\frac{c_2r}{v^{-1}(m/r)}\right\},\hspace{20pt}\forall m,r\in\mathbb{N}.\]}
\begin{proof} For any $m,r\in\mathbb{N}$, a standard argument (see proof of \cite{BarKum}, Theorem 4.9, for example) yields
\begin{equation}\label{twoterms}
p_{2m}^{T^*}(\rho,b_{2r})\leq\frac{P^{T^*}_\rho\left(X_{2m}=b_{2r},\tilde{\tau}_r\leq m\right)}{\mu^{T^*}(\{b_{2r}\})}+\frac{P^{T^*}_{b_{2r}}\left(X_{2m}={\rho},\tilde{\tau}_{r-1}\leq m\right)}{\mu^{T^*}(\{{\rho}\})},
\end{equation}
where $\tilde{\tau}_r$ is the stopping time for the random walk defined at (\ref{tildetau}). Applying the Markov property of $X$, we can bound the first of these terms as follows:
\begin{eqnarray*}
\lefteqn{\mu^{T^*}(\{b_{2r}\})^{-1}P^{T^*}_\rho\left(X_{2m}=b_{2r},\tilde{\tau}_r\leq m\right)}\\
&\leq &\mu^{T^*}(\{b_{2r}\})^{-1}E_\rho^{T^*}\left(\mathbf{1}_{\{\tilde{\tau}_r\leq m\}}P^{T^*}_{b_{r}}\left(X_{2m-\tilde{\tau}_r}=b_{2r}\right)\right)\\
&\leq & P^{T^*}_\rho\left(\tilde{\tau}_r\leq m\right) \sup_{m'\in[m,2m]}p_{m'}^{T^*}(b_r,b_{2r})\\
&\leq & P^{T^*}_\rho\left(\tilde{\tau}_r\leq m\right) \left(p_{2\lfloor m/2\rfloor}^{T^*}(b_r,b_{r})p_{2\lfloor m/2\rfloor}^{T^*}(b_{2r},b_{2r})\right)^{\tfrac{1}{2}},
\end{eqnarray*}
where we use the Cauchy-Schwarz inequality and the monotonicity of $p_{2k}(x,x)$ in $k$ to deduce the third inequality. Now, if $\gamma\in(0,1-\alpha^{-1})$, we can choose $\beta> 1$ large enough so that $\gamma(1+\beta)>1$ and also $\gamma':=\gamma(1+\beta)\beta^{-1}<1-\alpha^{-1}$. Consequently, applying H\"{o}lder's inequality with the exponents $1+\beta$, $2(1+\beta)/\beta$ and $2(1+\beta)/\beta$ to the random variables $P^{T^*}_\rho(\tilde{\tau}_r\leq m)^\gamma$, $p_{2\lfloor m/2\rfloor}^{T^*}(b_r,b_{r})^{\gamma/2}$ and $p_{2\lfloor m/2\rfloor}^{T^*}(b_{2r},b_{2r})^{\gamma/2}$ respectively, we have that
\begin{eqnarray}
\lefteqn{\mathbf{E}\left(\mu^{T^*}(\{b_{2r}\})^{-\gamma}P^{T^*}_\rho\left(X_{2m}=b_{2r},\tilde{\tau}_r\leq m\right)^{\gamma}\right)}\nonumber\\
&\leq&\mathbb{P}\left(\tilde{\tau}_r\leq m \right)^{\tfrac{1}{1+\beta}}\sup_{r'\in\mathbb{N}}\mathbf{E}\left(p_{2\lfloor m/2\rfloor}^{T^*}(b_{r'},b_{r'})^{\gamma'}\right)^{\tfrac{\beta}{1+\beta}},\label{tauexp}
\end{eqnarray}
where we also have used the fact that $P^{T^*}_\rho\left(\tilde{\tau}_r\leq m\right)^{\gamma(1+\beta)}\leq P^{T^*}_\rho\left(\tilde{\tau}_r\leq m\right)$. To bound the expectations in this expression, we proceed as at (\ref{hkupper}) to deduce the existence of a constant $c_1$ such that
\[\mathbf{E}\left(p_{2m}^{T^*}(b_{r},b_{r})^{\gamma'}\right)\leq \frac{c_1}{v(\mathcal{I}(m))^{\gamma'}}\mathbf{E}\left(1+\frac{v(R)^{\gamma'}}{{V(b_r,R)}^{\gamma'}}\right),\hspace{20pt}\forall m,r\in\mathbb{N},\]
where $V(b_r,R):=\{x\in T^*:d_{T^*}(b_r,x)\leq R\}$ and $R=R(m)$ is chosen to satisfy $\tfrac{1}{2}h(R)\leq m\leq h(R)$. By considering only the descendants of $b_r$, it is clear that $V(b_r,R)$ stochastically dominates $Y_R^*$ for every $r\in\mathbb{N}$. Thus, adjusting $c_1$ as necessary, it follows that
\[\sup_{r\in\mathbb{N}}\mathbf{E}\left(p_{2\lfloor m/2\rfloor}^{T^*}(b_{r},b_{r})^{\gamma'}\right)^{\tfrac{\beta}{1+\beta}}\leq\frac{c_1}{v(\mathcal{I}(m))^\gamma},\hspace{20pt}\forall m\in\mathbb{N}.\]
We now look to bound the first factor of the upper bound at (\ref{tauexp}). Recall the bound on the distribution of $\tilde{\tau}_r$ from (\ref{tail}) and let $c_2,c_3$ be the constants of this inequality. If $c_2h(r)\leq m$, then it is easy to check that $\exp(-r/v^{-1}(m/r))\geq c_4^{-1}>0$, thus \[\mathbb{P}(\tilde{\tau}_r\leq m) \leq c_4\exp(-r/v^{-1}(m/r))\]
in this case. We now assume that $c_2 h(r)>m$, choose $n$ to be the largest integer such that $c_2nh(\lfloor \tfrac{r}{n}\rfloor)\geq m$, and set $R=\lfloor \tfrac{r}{n}\rfloor$, so that (\ref{tail}) implies that
\[\mathbb{P}(\tilde{\tau}_r\leq m) \leq\mathbb{P}(\tilde{\tau}_{nR}\leq c_2nh(R)) \leq e^{-c_3 n}.\]
Applying this and the previous bound, it is elementary to check that there exist constants $c_5,c_6$ such that
\begin{equation}\label{exittime}
\mathbb{P}(\tilde{\tau}_r\leq m)^{\tfrac{1}{1+\beta}} \leq c_5\exp(-c_6 r/v^{-1}(m/r)),\hspace{20pt}\forall m,r\in\mathbb{N}.
\end{equation}
Thus we have so far demonstrated that
\[\mathbf{E}\left(\mu^{T^*}(\{b_{2r}\})^{-\gamma}P^{T^*}_\rho\left(X_{2m}=b_{2r},\tilde{\tau}_r\leq m\right)^{\gamma}\right)\leq \frac{c_7}{v(\mathcal{I}(m))^\gamma}\exp(-c_6 r/v^{-1}(m/r)),\]
for every $m,r\in\mathbb{N}$, for some finite constant $c_7$. To complete the proof, it remains to obtain a similar bound for the second term at (\ref{twoterms}), which can be done by following a similar argument to the one above. The one point that requires checking is that (\ref{exittime}) holds when $\mathbb{P}(\tilde{\tau}_r\leq m)$ is replaced by $\mathbb{P}_{b_{2r}}(\tilde{\tau}_{r-1}\leq m)$, where $\mathbb{P}_{b_{2r}}:=\int P^{T^*}_{b_{2r}}(\cdot)d\mathbf{P}$. Clearly, we have that $\mathbb{P}_{b_{2r}}(\tilde{\tau}_{r-1}\leq m)=\mathbb{P}_{b_{r+1}}(\tilde{\tau}_0\leq m)$, and so it will suffice to estimate the right-hand side of this inequality. Now define, for each $r\in\mathbb{N}$, the subset
\[T^*_r:=\{x\in T^*:x\mbox{ is $b_r$ or a descendant of $b_r$}\},\]
and set $T^*(r)=T^*\backslash T^*_{r+1}$. By the description of $T^*$ in \cite{Kesten}, Lemma 2.2, we know that the subtrees growing out of the backbone of $T^*$ are independent and identically distributed. Applying this fact, it is an elementary exercise to check that the law of $\tilde{\tau}_0$ under $\int P^{T^*(r)}_{b_r}(\cdot)d\mathbf{P}$ is the same as the law of $\tilde{\tau}_r$ under $\int P^{T^*(r)}_{\rho}(\cdot)d\mathbf{P}$, which in turn is the same as the law of $\tilde{\tau}_r$ under $\mathbb{P}$. In particular, we have that $\int P^{T^*(r)}_{b_r}(\tilde{\tau}_0\leq m)d\mathbf{P}=\mathbb{P}(\tilde{\tau}_r\leq m)$ for every $m,r\in\mathbb{N}$. Furthermore, by construction, we have that the left-hand side of this identity is equal to $\mathbb{P}_{b_r}(\tilde{\tau}_0'\leq m)$, where
\[\tilde{\tau}_0':=\#\{n\in[0,\tilde{\tau}_0):X_n,X_{n+1}\in T^*(r)\}.\]
Since $\tilde{\tau}_0\geq \tilde{\tau}_0'$, it follows that $\mathbb{P}_{b_r}(\tilde{\tau}_0\leq m)\leq\mathbb{P}_{b_r}(\tilde{\tau}_0'\leq m)=\mathbb{P}(\tilde{\tau}_r\leq m)$, and the result follows.
\end{proof}

In the case when the offspring distribution is binomial, it is straightforward to check that the upper bound deduced in the above proposition is sharp up to constants by applying estimates of \cite{BarKum}. In general, however, we are only able to prove the corresponding lower bound holds near the diagonal. That we can not extend the chaining argument of \cite{BarKum} to obtain the full off-diagonal lower bound (even along the backbone) results from the fact that we only have a polynomial tail bound for the probability that $T^*$ admits ``bad'' subsets, whereas, in the binomial case, proving an exponential tail bound is possible.

{\propn If $\gamma>0$, then there exist constants $c_1,c_2\in(0,\infty)$ such that
\[\mathbf{E}\left(p_{2m}^{T^*}(\rho,b_{2r})^\gamma\right)^{\tfrac{1}{\gamma}}\geq c_1 v(\mathcal{I}(m))^{-1},\]
whenever $1\leq r\leq c_2\mathcal{I}(m)$.}
\begin{proof} By a standard argument (cf. \cite{BarKum}, Proposition 4.4), there exists a deterministic constant $c_1$ such that if $T^*$ satisfies, for some $R\geq 2$, $\lambda\geq 8$,
\begin{equation}\label{conds1}
V(\lambda R)\in[\lambda^{-1}v(\lambda R),\lambda v(\lambda R)],\hspace{10pt}V(R)\geq \lambda^{-1}v(R),\hspace{10pt}R_{eff}(B(R),B(\lambda R)^c)\geq 4R,
\end{equation}
then, for $m\in [\tfrac{1}{2}\lambda^{-1}h(R-1),\tfrac{1}{2} \lambda^{-1}h(R)]$,
\[p^{T^*}_{2m}(x,x)\geq c_1\lambda^{-\theta_1}v(\mathcal{I}(m))^{-1},\hspace{20pt}\forall x\in B(R),\]
where $\theta_1:=19/(\alpha-1)$. This is easily extended (cf. \cite{BarKum}, Theorem 4.6(c)) to the result that
\[p^{T^*}_{2m}(\rho,b_{2r})\geq c_2\lambda^{-\theta_1}v(\mathcal{I}(m))^{-1},\]
for every $1\leq r\leq c_3\lambda^{-\theta_1}\mathcal{I}(m)$, for some constants $c_2$ and $c_3$. A straightforward adaptation of the proof of Lemma \ref{ass12} allows it to be proved that the conditions at (\ref{conds1}) hold with probability greater than $\tfrac{1}{2}$ for some $\lambda\geq \max\{8,v(1)\}$, uniformly in $R\geq 2$. Using this choice of $\lambda$, if $m\in\mathbb{N}$, we can choose $R\geq 2$ that satisfies $m\in [\tfrac{1}{2}\lambda^{-1}h(R-1),\tfrac{1}{2} \lambda^{-1}h(R)]$, and applying the lower bound above, it follows that
\[\mathbf{E}\left(p_{2m}^{T^*}(\rho,b_{2r})^\gamma\right)^{\tfrac{1}{\gamma}}\geq c_4 v(\mathcal{I}(m))^{-1},\]
for $1\leq r\leq c_5\mathcal{I}(m)$, which completes the proof.
\end{proof}

Summarising the two previous results using the expressions for $v$ and $\mathcal{I}$ presented in the introduction, we have the following bounds for the transition density of the simple random walk on $T^*$.

{\cor If $\gamma\in(0,1-\alpha^{-1})$, then there exist constants $c_1,c_2\in(0,\infty)$ and slowly varying functions $\ell_1$, $\ell_2$ and $\ell_3$ such that
\[\mathbf{E}\left(p_{2m}^{T^*}(\rho,b_{2r})^\gamma\right)^{\tfrac{1}{\gamma}}\leq c_1 m^{-\tfrac{\alpha}{2\alpha-1}}\ell_1(m) \exp\left\{-\left(\frac{r^{\tfrac{2\alpha-1}{\alpha-1}}}{m}\right)^{\tfrac{\alpha-1}{\alpha}}\ell_2\left(\frac{m}{r}\right)\right\},\]
for every $m,r\in\mathbb{N}$, and also
\[\mathbf{E}\left(p_{2m}^{T^*}(\rho,b_{2r})^\gamma\right)^{\tfrac{1}{\gamma}}\geq c_2 m^{-\tfrac{\alpha}{2\alpha-1}}\ell_1(m),\]
whenever $1\leq r\leq m^{-\tfrac{\alpha}{2\alpha-1}}\ell_3(m)$.}

\section{Volume and transition density fluctuations}\label{fluct}

To establish that the transition density of $X$ exhibits logarithmic fluctuations when $\alpha\in(1,2)$, and at least log-logarithmic fluctuations when $\alpha=2$, which is the aim of this section, we start by showing that the same is true of the volume growth on the tree $T^*$.

{\lem (a) If $\beta_1\in(0,\alpha-1)$, then $\mathbf{P}$-a.s. realisation of $T^*$ satisfies
\[\limsup_{R\rightarrow\infty}\frac{V(R)}{v(R)(\log R)^{1/\beta_1}}=0.\]
(b) If $\alpha\in(1,2)$ and $\beta_2>(\alpha-1)/(2-\alpha)$, then $\mathbf{P}$-a.s. realisation of $T^*$ satisfies
\[\limsup_{R\rightarrow\infty}\frac{V(R)}{v(R)(\log R)^{1/\beta_2}}=\infty.\]
If $\alpha=2$ and $\varepsilon>0$, then $\mathbf{P}$-a.s. realisation of $T^*$ satisfies
\[\limsup_{R\rightarrow\infty}\frac{V(R)}{v(R)(\log\log R)^{1-\varepsilon}}=\infty.\]}
\begin{proof} Clearly, by (\ref{215}), it will suffice to prove the result with $Y^*_R$ in place of $V(R)$. By the Borel-Cantelli lemma, part (a) is an easy consequence of Proposition \ref{yupper}. To prove (b), we consider the sequence of subsets $(A_n)_{n\geq 0}$ of $T^*$ defined by
\[A_n:=\{x\in T^*:d_{T^*}(\rho,x)\in[2^n,2^{n+1}),\mbox{ $x$ is $b_{2^n}$ or a descendant of $b_{2^n}$}\},\]
where $b_{2^n}$ is the point on the backbone at a distance $2^n$ from the root. By \cite{Kesten}, Lemma 2.2, $(\#A_n)_{n\geq 0}$ is a sequence of independent random variables, and $\#A_n$ is equal in distribution to $Y_{2^{n+1}-2^n-1}^*$. Thus, if $\alpha\in(1,2)$ and $\beta_2>(\alpha-1)/(2-\alpha)$, we can apply the second Borel-Cantelli lemma and (\ref{yupperc}) to obtain that $Y^*_{2^{n+3}}\geq \#A_{n+2}\geq v(2^n)n^{1/\beta_2}$ infinitely often, $\mathbf{P}$-a.s., and the first claim follows. In the case when $\alpha=2$, note that Proposition \ref{zlower} allows us to choose strictly positive constants $c_1$ and $c_2$ such that $\mathbf{P}(Z_n^*>c_1p_n^{-1})\geq c_2$ for every $n\in\mathbb{N}$. By considering a decomposition of $T^*$ similar to the one applied in the proof of Proposition \ref{yupper}, it follows that, for every $n\in\mathbb{N}$, $\lambda>0$,
\begin{eqnarray*}
\mathbf{P}\left(Y^*_{2n}\geq \lambda np_n^{-1}\right)&\geq& c_2\mathbf{P}\left(Y^*_{2n}\geq \lambda np_n^{-1}|Z_n^*>c_1p_n^{-1}\right)\\
&\geq &c_2\mathbf{P}\left(\mathrm{Bin}(\lfloor c_1p_n^{-1}\rfloor,c_3p_n)\geq c_4\lambda \right),
\end{eqnarray*}
for suitably chosen $c_3,c_4$. Straightforward estimates (cf. \cite{BarKum}, (2.18)) subsequently imply that
\[\mathbf{P}\left(Y^*_{2^{n+1}}\geq v(2^n)(\log n)^{1-\varepsilon}\right)\geq c_5n^{-1},\hspace{20pt}\forall n\in\mathbb{N},\]
for some $c_5>0$. This estimate allows us to apply the second Borel-Cantelli lemma, as in the case $\alpha\in(1,2)$, to complete the proof.
\end{proof}

In addition to the above lemma, note that Fatou's Lemma and the moment estimate for $Y_R^*$ at (\ref{ymoments}) implies that $\liminf_{R\rightarrow\infty}V(R)/v(r)<\infty$, $\mathbf{P}$-a.s. Hence there are $\mathbf{P}$-a.s. asymptotic fluctuations about $v(R)$ in the volume growth on $T^*$ of at least log-logarithmic order when $\alpha=2$ and of at least logarithmic order when $\alpha\in(1,2)$. Furthermore, from the previous lemma we are immediately able to determine the following asymptotic result for the transition density of $X$, which can be proved in the same way as \cite{BarKum}, Lemma 5.1. In conjunction with (\ref{anntd}), these results demonstrate that with positive probability there are fluctuations in the transition density about $v(\mathcal{I}(m))^{-1}$.

{\cor  If $\alpha\in(1,2)$, then there exists an $\varepsilon_1>0$ such that, $\mathbf{P}$-a.s., the transition density of $X$ satisfies
\[\liminf_{m\rightarrow\infty} v(\mathcal{I}(m))(\log m)^{\varepsilon_1}p^{T^*}_{2m}(\rho,\rho)=0.\]
If $\alpha=2$, then there exists an $\varepsilon_2>0$ such that, $\mathbf{P}$-a.s., the transition density of $X$ satisfies
\[\liminf_{m\rightarrow\infty} v(\mathcal{I}(m))(\log\log m)^{\varepsilon_2}p^{T^*}_{2m}(\rho,\rho)=0.\]}

%TK new
\section*{Acknowledgements}
The authors thank the anonymous referees for carefully reading the first draft and
for suggesting various improvements.
%end new

\def\cprime{$'$}
\providecommand{\bysame}{\leavevmode\hbox to3em{\hrulefill}\thinspace}
\providecommand{\MR}{\relax\ifhmode\unskip\space\fi MR }
\providecommand{\MRhref}[2]{%
  \href{http://www.ams.org/mathscinet-getitem?mr=#1}{#2}
}
\providecommand{\href}[2]{#2}

\end{document}